\documentclass[12pt]{amsart}
\usepackage{amssymb}

\usepackage{euscript}

\let\cal\mathcal
\input xypic
\xyoption {all}

\makeatletter \@mparswitchfalse \makeatother

\newtheorem{theorem}{Theorem}
\newtheorem{lemma}[theorem]{Lemma}
\newtheorem{corollary}[theorem]{Corollary}
\newtheorem{proposition}[theorem]{Proposition}
\newtheorem{remarks}[theorem]{Remarks}
\newtheorem{sublemma}[theorem]{Sublemma}
\newtheorem{definition}[theorem]{Definition}

\def\m{{\cal M}}
\def\M{{\cal M}}

\def\mt{\widetilde{\cal M}}
\def\p{\prec\prec}
\def\ch{\raise 0.5ex \hbox{$\chi$}}

\def\nmt{\widetilde {\cal M}}

\def\nm{{\cal M}}
\def\N{{\rm I\kern-.1567em N}}
\def\3{\ss}
\def\R{{\rm I\kern-.1567em R}}

\def\mp+{$L_p(\nm)+{\cal M}$}
\def\mp+n{L_p(\nm)+{\cal M}}

\numberwithin{equation}{section}
\numberwithin{theorem}{section}

\title[non-commutative spaces]{Embeddings of $\ell_{\boldsymbol p}$
into non-commutative spaces}

\author{Narcisse Randrianantoanina}
\address{Department of Mathematics and Statistics, Miami University, Oxford,
Ohio 45056}
\thanks{Supported in part by NSF Grant DMS-9703789}
\email{randrin@muohio.edu}

\begin{document}

 \begin{abstract} Let $\M$ be a semi-finite von Neumann algebra
equipped with a faithful normal trace $\tau$. We study the 
subspace structures  of non-commutative Lorentz spaces $L_{p,q}(\M, \tau)$,
extending  results of Carothers and Dilworth to the non-commutative settings.
In particular, we show that,
under natural conditions on indices, $\ell_p$ can not be
embedded into $L_{p,q}(\M, \tau)$. As applications, we prove that 
for $0<p<\infty$ with $p \neq 2$ then $\ell_p$ cannot be strongly embedded into $L_p(\M,\tau)$.
Thus providing a non-commutative extension of 
a result of Kalton for $0<p<1$ and a result of
Rosenthal for $1\leq p <2$  on $L_p[0,1]$.
\end{abstract}

\maketitle

\section{Introduction}
 The study of rearrangement invariant
 Banach spaces of measurable functions is a classical theme.
Several studies have been devoted on characterizations of subspaces of 
rearrangement invariant spaces. 
Recently,  the theory of rearrangement invariant Banach spaces
of measurable operators affiliated with semi-finite von Neumann algebra
have emerged as the natural non-commutative generalizations of 
K\"othe functions spaces. This theory, which is based on
the theory of non-commutative integration introduced by Segal
(see \cite{SEG}), replaces the classical duality $(L_\infty(\mu), L_1(\mu))$
by the duality between a semi-finite von Neumann algebra and its predual.
 It  provides a unified approach  to the study of unitary ideals and
 rearrangement invariant spaces. Several authors have considered 
these non-commutative spaces of measurable operators (see for instance,
\cite{DDP1}, \cite{DDP3}, \cite{D3LS}, \cite{CS} and \cite{KOS2}).

The purpose of the present paper is to examine the subspaces
of symmetric spaces of measurable operators in which the norm topology and 
the measure topology coincide and subspaces generated by disjointly
supported basic sequences.
The classical spaces $L_p(\mu)$ are of central
inportance and  
results in their  structures go back to  the work of Banach.
Since their introduction by Lorentz in 1950, the Lorentz function 
spaces $L_{p,q}$ have been found to be of special interests in many aspects
of analysis and probability theory. In \cite{CADI2} and \cite{CADI}, 
Carothers and Dilworth
studied the spaces $L_{p,q}[0,1]$ and $L_{p,q}[0,\infty)$. They proved, 
among other things, that for some appropriate values of the indices 
$p$ and $q$,  $L_{p,q}[0,\infty)$ does not contain $\ell_p$. Presisely:
\begin{theorem}
Let $0<p<\infty$, $0<q<\infty$, $p\neq q$, $p\neq 2$. 
Then $\ell_p$ does not embed
into $L_{p,q}[0,\infty)$.
\end{theorem}
Motivated 
by  Theorem~1.1, we examine the subspace
 structure of 
the  non commutative Lorentz spaces $L_{p,q}(\M,\tau)$, where $(\M,\tau)$ is 
a semi-finite von 
Neumann algebra. The principal result of the present paper is that the main
results of \cite{CADI2} and \cite{CADI} extend  to the non-commutative 
settings.
 The initial basic question, that led to the consideration of these
Lorentz spaces, is the question of embeddings of $\ell_p$ into 
$L_p(\cal M, \tau)$. Clearly, any disjointly supported basic sequence in
$L_p(\cal M, \tau)$ is isomorphic to $\ell_p$. For the commutative case,
Rosenthal proved in \cite{R6} that  if $(\Omega,\Sigma,\mu)$ is a 
$\sigma$-finite measure space, $1\leq p <2$ and $X$ be a subspace of 
$L_p(\Omega,\Sigma,\mu)$ containing $\ell_p$ then the norm topology
and the measure topology do not coincide on $X$. For $0<p<1$, 
the same result can be
found implicitely in a paper of Kalton \cite{KA3}.
This implies that for $0<p<2$, any basic sequence in 
$L_p(\Omega,\Sigma,\mu)$ that is equivalent to $\ell_p$
is essentially a perturbation of a disjointly supported  basic sequence. 
We establish, as applications of our results on Lorentz spaces, that 
Kalton and Rosenthal's results extend to  $L_p(\m, \tau)$.

The paper is organized as follows. In Section~2 below, we gather some 
necessary definitions and present some basic facts conserning symmetric 
spaces of measurable operators that will be needed throughout.
In particular, we prove a Kadec-Pe\l czy\'nski type dichotomy for 
basic sequences in symmetric
spaces of measurable operators  of Rademacher type~2, generalizing a
result of Sukochev (see \cite{SU}). The final section is devoted entirely
to the study of subspaces of Lorentz spaces and its applications. 
In particular, Theorem~1.1 is generalized to the non-commutative case.
Our approach relies on a disjointification techniques based on the 
non-commutative analogue of the Khintchine's inequalities  of
Lust-Piquard and Pisier (\cite{LP4} and \cite{LPI}).

\section{Definitions and  Preliminaries}
We begin by recalling some definitions and facts about function spaces.
Let $E$ be a complex quasi-Banach lattice. If $0<p<\infty$, then $E$
is said to be {\it $p$-convex} (respectively {\it $p$-concave}) if 
there exists a constant $C>0$ such that all finite sequence 
$\{x_n\}$ in $E$,
$$\left\Vert \left(\sum \vert x_n \vert^p \right)^{\frac{1}{p}} 
\right\Vert_E
\leq C \left( \sum \left\Vert x_n \right\Vert^p \right)^{\frac{1}{p}}$$
$$\left( \text{resp.} \ 
\left\Vert \left(\sum \vert x_n \vert^p \right)^{\frac{1}{p}} 
\right\Vert_E
\geq C^{-1} \left( \sum \left\Vert x_n \right\Vert^p \right)^{\frac{1}{p}}
\right).$$
The least constant $C$ is called the $p$-convexity (respectively 
$p$-concavity) constant of $E$ and is denoted by $M^{(p)}(E)$ 
(respectively $M_{(p)}(E)$).

For $0<p<\infty$, $E^{(p)}$ will denote the quasi-Banach lattice defined by
$$E^{(p)} :=\left\{x\ :\  \vert x \vert^p \in E \right\}$$
equipped with 
 $$\left\Vert x \right\Vert_{E^{(p)}} = \left\Vert  \vert x \vert^p
\right\Vert_{E}^{1/p}.$$
It is easy to verify that if $E$ is $\alpha$-convex and $q$-concave then 
$E^{(p)}$ is $\alpha p$-convex and $qp$-concave with 
$M^{(\alpha p)}( E^{(p)})\leq M^{(\alpha)}(E)^{1/p}$ and 
$M_{(q p)}( E^{(p)})\leq M_{(q)}(E)^{1/p}$. Consequently, if $E$ is 
$\alpha$-convex then $E^{(1/\alpha)}$ is $1$-convex and therefore can be 
equivalently renormed to be a Banach lattice (\cite{LT}).

The quasi-Banach lattice $E$ is said to satisfy {\it a lower $q$-estimate}
(respectively {\it upper $p$-estimate}) if there exists a positive constant
$C>0$ such that for all finite sequence of mutually disjoint elements of $E$
$$(\sum \Vert x_n \Vert_{E}^q)^{1/q} \leq C \Vert \sum x_n \Vert_E$$
$$\left( \text{resp.}\ 
(\sum \Vert x_n \Vert_{E}^p)^{1/p} \geq C^{-1} \Vert \sum x_n \Vert_E
\right).$$
 We denote by
${\cal M}$ a semi-finite von Neumann algebra on the Hilbert space
${\cal H}$, with a fixed faithful and normal semi-finite trace $\tau $. The
identity in ${\cal M}$  is denoted by ${\bf 1}$, and we denote by
$ {\cal M}^p$ the
set of all projections in ${\cal M}$. A linear operator
$x$:dom$(x)\to {\cal H} $, with domain
dom$(x)\subseteq {\cal H}$, is called {\it affiliated with}  $\m$  if
$ux=xu$ for all unitary $u$ in the commutant $\nm'$ of $\m.$
The closed and densely defined operator $x$, affiliated with
${\cal M}$ , is called
$\tau$-{\it measurable} if for every $\epsilon >0$ there exists
$p\in \nm^p$ such the $p({\cal H})\subseteq$dom$(x)$ and
$\tau({\bf 1}-p)<\epsilon$.  With the sum and product defined as
the respective
closures of the algebraic sum and product, $\mt $ is a *-algebra.
 For standard facts concerning von Neumann
algebras, we refer to \cite{KR} and \cite{TAK}.
\smallskip
\noindent
We recall the notion of generalized singular value function \cite{FK}.
Given a self-adjoint operator $x$ in ${\cal H}$ we denote by $e^x(\cdot
)$ the spectral measure of $x$. Now assume that $x\in\nmt$.
 Then $e^{\vert x\vert}(B)\in\nm$ for all Borel sets $B\subseteq{\mathbb R}$,
and there exists $s>0$ such that
$\tau (e^{\vert x\vert}(s,\infty))<\infty$.
 For $x\in\nmt$ and $t\geq 0$ we define
$$
\mu_t(x)=\inf\{s\geq 0 : \tau (e^{\vert x\vert}(s,\infty))\leq t\}.
$$
The function $\mu(x):[0,\infty)\to [0,\infty ]$ is called the {\it
generalized singular value function} (or decreasing rearrangement) of
$x$; note that $\mu_t(x)<\infty$ for all $t>0$.
 Suppose that $a >0$. If we consider $\nm =L_\infty([0,a),m)$,
 where $m$ denotes Lebesgue
measure on the interval $[0,a)$, as an Abelian von Neumann algebra
acting via
multiplication on the Hilbert space  ${\cal H}=L_2([0,a),m)$,
with the
trace given by integration with respect to $m$, it is easy to see that
$\nmt $ consists of all measurable functions on $[0,a )$ which are
bounded
except on a set of finite measure. Further, if $f\in\nmt$, then the
generalized singular value function $\mu(f)$ is precisely the classical
non-increasing rearrangement of the function $\vert f\vert$.
On the other hand, if $({\cal M},\tau)$ is the space of all
bounded linear operators in some Hilbert space equipped with
the canonical trace $tr$, then $\nmt =\nm$ and, if $x\in \nm$
is compact, then the generalized singular value function $\mu (x)$
may be identified in a natural manner with the sequence
$\{\mu _n(x)\}_{n=0}^\infty$ of singular values of
$\vert x\vert =\sqrt {x^*x}$, repeated according to
multiplicity and arranged in non-increasing order.
\smallskip
\noindent  By $L_0([0, a),m)$ we denote the
space of all ${\mathbb C}$-valued Lebesgue measurable functions on the
interval $[0,a )$
(with identification $m$-a.e.).
 A quasi-Banach space $(E,\Vert \cdot
\Vert_{_E})$, where $E\subseteq L_0([0,a),m)$ is called a {\it
rearrangement-invariant Banach function space} on the interval
$[0,a)$, if it follows from $f\in
E, g\in L_0([0,a),m)$ and $\mu(g)\leq\mu (f)$ that $g\in E$
and $\Vert g\Vert_E\leq \Vert f\Vert_E$. If $(E,\Vert\cdot\Vert_E)$
is a rearrangement-invariant  quasi-Banach function space on $[0,a)$,
then $E$ is said to be {\it symmetric }
if $f,g\in E$ and $g\p f$ imply
that $\Vert g\Vert_E\leq\Vert f\Vert_E$.
 Here $g\p f$ denotes submajorization in the sense of
 Hardy-Littlewood-Polya :
$$
\int^t_0 \mu_s(g) ds \leq  \int^t_0 \mu_s(f) ds,
\quad {\rm for \ all}\ t>0.
$$
%\noindent If $(E,\Vert\cdot\Vert_E)$ is a rearrangement-invariant, symmetric
%Banach function space on $[0,\alpha )$, then $E$ will
%be called {\it fully symmetric} if
%$f\in E, g\in L_0([0,\alpha),m)$ and $g\prec\prec f$ implies $g\in E$
%and $\Vert g\Vert _{_{E}}\leq \Vert f\Vert _{_{E}}$.
The general theory of rearrangement-invariant spaces may
be found in \cite{LT} and \cite{SA}.

\noindent Given a semi-finite von Neumann algebra $(\nm,\tau)$
and a symmetric quasi-Banach function space
$(E,\Vert \cdot \Vert _{_{E}})$ on $([0,\tau ({\bf 1})),m)$,
we define the non-commutative space $E(\nm,\tau)$ by setting
$$E(\nm,\tau):=\{x\in \nmt: \mu (x)\in E\}\quad \text{with}$$
$$\Vert x\Vert _{_{E(\nm,\tau)}}:=\Vert \mu(x)\Vert _{_{E}}\ 
\text{for}\ x \in E(\nm,\tau).$$
It is known that if $E$ is $\alpha$-convex  for some $0<\alpha<\infty$
with $M^{(\alpha)}(E)=1$ then $\Vert \cdot \Vert _{_{E(\m,\tau)}}$
is a norm for $\alpha \geq 1$ and an $\alpha$-norm  if
$0<\alpha<1$. 
 The space 
$(E(\nm,\tau), \Vert \cdot \Vert _{_{E(\nm,\tau)}})$ is a $\alpha$-Banach
space.
 Moreover,  the inclusions
$$
 L_\alpha({\cal M}, \tau ) \cap {\cal M}
 \subseteq E ({\cal M}, \tau )
 \subseteq L_\alpha({\cal M}, \tau ) + {\cal M}.
$$
hold with continuous embeddings. 
 We remark that if 
$0<p<\infty$ and 
$E=L_p([0,\tau({\bf 1}))$ then $E(\cal M, \tau)$ coincides with the definition
of $L_p(\cal M,\tau)$ as in \cite{N} and \cite{Te}. In particular, if
$\m =\cal{L}(\cal H)$ with the standard trace then these
 $L_p$-spaces are precisely the {\it Schatten classes} $\cal C_p$.

We recall that 
the topology defined by the metric on $\widetilde{\m}$ obtained by setting
$$d(x,y)=\inf\left\{t\geq 0:\ \mu_t(x-y)\leq t \right\}, 
\ \ \ \text{for}\ 
x, y \in \widetilde{\m},$$
is called the {\it measure topology}. It is well-known that a net
$(x_\alpha)_{\alpha \in I}$ in $\widetilde{\m}$  converge to 
$x \in \widetilde{\m}$ in measure  topology if and only if
for every $\epsilon >0$, $\delta >0$, there exists $\alpha_0 \in I$
such that whenever $\alpha \geq \alpha_0$, there exists a projection
$p \in \m^p$ such that
   $$\left\Vert (x_\alpha - x)p \right\Vert_{\m} <\epsilon \ \
\text{and}\ \ \tau({\bf 1}-p)<\delta.$$
It was shown in \cite{N} that $(\widetilde{\m}, d)$ is 
a complete metric, Hausdorff, topological ${}^*$-algebra.
For $x \in \widetilde{\m}$, the right and left support projections of $x$
are denoted by $r(x)$ and $l(x)$ respectively. Operators 
$x, y \in \widetilde{\m}$ are said to be right (respectively, left)
disjointly supported if $r(x)r(y)=0$ (respectively, $l(x)l(y)=0$).
\smallskip

\begin{definition} 
Let $E$ be a symmetric quasi-Banach function space on $[0,\tau({\bf 1}))$.
We say that a subspace $X$ of $E(\m,\tau)$ is strongly embedded into 
$E(\m,\tau)$
if the $\Vert \cdot \Vert_{E(\m,\tau)}$-topology and the measure topology
on $X$ coincide.
\end{definition}

\smallskip

The following definition was introduced in \cite{Ran10} as an analogue of
the uniform integrability of family of functions.
\begin{definition}
Let  $E$ be a symmetric  quasi-Banach function space on $[0,\tau({\bf 1}))$.
A bounded subset $K$ of $E(\cal{M}, \tau)$ is said to be
{ $E$-uniformly-integrable} if $\lim\limits_{n \to \infty} \sup\limits_{x \in K}
\left\Vert e_n x e_n \right\Vert_{E(\cal{M},\tau)} = 0$
for every decreasing sequence $\{e_n\}_{n=1}^\infty$ of projections with
$e_n \downarrow_n  0$.
\end{definition}

A non-commutative extension of the Kadec-Pe\l czy\'nski subsequence 
decomposition relative to the above notion of uniform integrability
was considered in \cite{Ran10} (see Theorem~3.1, Theorem~3.7 and
Corollary~3.8) and will be used repeatedly throughout 
this paper.

\smallskip

For the remaining of this section, we will present some results, 
some of which may be of independent interests, that we will need 
in the sequel.

\smallskip

The following proposition is essentially due to Sukochev
(\cite{SU}).
\begin{proposition}\label{disjoint}
Let $E$ be $\alpha$-convex with constant $1$ and assume that
$E$ is order continuous.
Let $\{x_{n}\}_{n=1}^\infty$ be a basic sequence in $E(\m, \tau)$ such that
$\{x_{n}\}_{n=1}^\infty$ is both right and left disjointly supported then
$\{x_{n}\}_{n=1}^\infty$ is equivalent to a disjointly supported
 basic sequence in $E$.
\end{proposition}
\begin{proof} For each $n\geq 1$, 
let $q_{n} := l(x_{n})$ and $p_{n} := r(x_{n})$ be the left and right
support projection of $x_{n}$ respectively. Both sequences
$\{q_{n}\}_{n=1}^\infty$ and $\{p_{n}\}_{n=1}^\infty$ are mutually disjoint  and for every
$n\geq 1$,  $x_{n} = q_{n} x_{n}
p_{n}$.  For any finite sequence of scalars $\{a_i\}_{i=1}^n$,
\begin{equation*}
\begin{split}
\left\vert \sum^{n}_{i=1} a_{i} x_{i} \right\vert^{2} &=\left( \sum^{n}_{i=1}
\bar{a}_{i} p_{i} x_{i}^{*} q_{i} \right) \left( \sum^{n}_{i=1} a_{i} q_{i}
x_{i} p_{i} \right) \cr
&= \sum^{n}_{i=1} \vert a_{i} \vert^{2} p_{i} x_{i}^{\ast} q_{i} x_{i}
p_{i} \cr
&= \left\vert \sum_{i=1}^n a_{i} \vert x_{i} \vert \right\vert^{2}.
\end{split}
\end{equation*}
Note that $\{\vert x_{i} \vert\}_{i=1}^\infty$ is disjointly supported by the
projections $\{p_{i}\}_{i=1}^\infty$. For each $i\geq 1$,
the semi-finiteness of $p_{i}$ implies that the  family $\{e_\beta\}_\beta$ of all
projections in $p_i\m p_i$ of finite trace  satisfies
$0\leq e_\beta \uparrow_\beta p_i$. Since $E$ is order-continuous, it
follows that 
$\left\Vert e_\beta \vert x_i \vert e_\beta 
- \vert x_i \vert \right\Vert \rightarrow_\beta 0$.

For each $i \geq 1$, choose a projection $\tilde{p}_{i} \leq p_{i}$ such
that $\tau(\tilde{p}_{i}) < \infty$ and 
$\Vert \tilde{p}_{i} \vert
x_{i} \vert \tilde{p}_{i} - \vert x_{i} \vert \Vert^\alpha \leq 2^{-i}$.

\noindent 
{\it Claim:\ The sequence 
 $\{\tilde{p}_{i} \vert x_{i} \vert \tilde{p}_{i}\}_{i=1}^\infty$ is
equivalent to $\{\vert x_{i} \vert\}_{i=1}^\infty$.} 

Let $p= \vee_{i=1}^\infty \tilde{p}_i$. For any 
$x=\sum_{i=1}^\infty a_i\vert x_i\vert \in 
\overline{\text{span}}\{\vert x_i \vert, i\geq 1\}$, we have
$\sum_{i=1}^\infty a_i \tilde{p}_{i} \vert x_{i} \vert \tilde{p}_{i}=
p(\sum_{i=1}^\infty a_i\vert x_i\vert)p$ so the series 
$ \sum_{i=1}^\infty a_i \tilde{p}_{i} \vert x_{i} \vert \tilde{p}_{i}$
is convergent whenever $\sum_{i=1}^\infty a_i\vert x_i\vert$ does.
Conversely, if $\{a_n\}_{n=1}^\infty$ is a bounded sequence of
scalars such that 
$\sum_{i=1}^\infty a_i \tilde{p}_{i} \vert x_{i} \vert \tilde{p}_{i}$
is convergent, then for any subset $S$ of $\mathbb N$,
\begin{equation*}
\begin{split}
\left\Vert \sum_{i \in S} a_i \vert x_i \vert 
\right\Vert_{E(\m,\tau)}^\alpha 
&\leq \left\Vert \sum_{i \in S} a_i (\tilde{p}_i\vert x_i \vert 
\tilde{p}_i)\right\Vert_{E(\m,\tau)}^\alpha + 
\left\Vert \sum_{i \in S} a_i (\tilde{p}_i\vert x_i \vert 
\tilde{p}_i - \vert x_i \vert)\right\Vert_{E(\m,\tau)}^\alpha \cr
&\leq \sup_{i\in S}\vert a_i \vert^\alpha \cdot \sum_{i\in S} 2^{-i}
+ \left\Vert \sum_{i\in S} a_i \tilde{p}_{i} \vert x_{i} \vert 
\tilde{p}_{i} \right\Vert_{E(\m,\tau)}^\alpha.
\end{split}
\end{equation*}
This shows that the series $\sum_{i=1}^\infty a_i \vert x_i \vert$ is 
convergent.
Let $C_1$ and $C_2$ be positive constants so that  for any 
finite sequence of scalars
$\{a_i\}_{i=1}^n$,
$$C_1  \left\Vert\sum^{n}_{i=1} a_{i}\vert x_{i}\vert 
\right\Vert_{E(\m, \tau)} \leq
\left\Vert \sum^{n}_{i=1}a_{i} \tilde p_{i} \vert x_i \vert \tilde
p_{i} \right\Vert_{E(\m, \tau)} \leq
C_2 \left\Vert\sum^{n}_{i=1} a_{i}\vert x_{i}\vert \right\Vert_{E(\m, \tau)}.$$ 
If $\alpha_{1}=0$ and $\alpha_{n} = \sum \limits_{i=1}^{n} \tau(\tilde{p}_{i})
< \infty$,  set $f_{n}: = \mu_{(\cdot)- \alpha_{n}} (\tilde{p}_{n} \vert
x_{n} \vert \tilde{p}_{n})$ for $n \geq1$.  The sequence 
$\{f_{n}\}_{n=1}^\infty$
is disjointly supported in $E(0, \tau(\bf{1}))$ and 
$\{f_{n}\}_{n=1}^\infty$ is
isometrically isomorphic to $\{\tilde{p}_{n} \vert x_{n} \vert
\tilde{p}_{n}\}_{n=1}^\infty$. For any finite sequence of scalars $\{a_i\}_{i=1}^n$,
\begin{equation*}
\begin{split}
C_1\left\Vert \sum^{n}_{i=1} a_{i} x_{i} \right\Vert_{E(\m, \tau)} &= C_1
\left\Vert
\sum^{n}_{i=1} a_{i}\vert x_{i}\vert \right\Vert_{E(\m, \tau)} \cr
&\leq \left\Vert \sum^{n}_{i=1}a_{i} \tilde p_{i} \vert x_i \vert \tilde
p_{i} \right\Vert_{E(\m, \tau)} \cr
&= \left\Vert \sum^{n}_{i=1} a_{i} f_{i} \right\Vert_{E} \cr
&\leq C_2\left\Vert \sum^{n}_{i=1} a_{i} x_{i} \right\Vert_{E(\m, \tau)} .
\end{split}
\end{equation*}
The proof is complete.
\end{proof}

\begin{proposition}\label{dichotomy1}
Let $E$ be an order continuous symmetric quasi-Banach function space on
$[0,\tau({\bf 1}))$ that is $\alpha$-convex with constant $1$ for some 
$0<\alpha\leq 1$ and suppose that $E$ satisfies a lower $q$-estimate 
with constant $1$ for some $q\geq \alpha$. If $X$ is a subspace of 
$E(\m,\tau)$ then either $X$ is strongly embedded into 
$E(\m,\tau)$ or there exists a normalized basic sequence 
$\{y_n\}_{n=1}^\infty$ in $X$,
a mutually disjoint sequence of projections $\{e_n\}_{n=1}^\infty$ in $\m$ such that
$$\lim_{n\to \infty} \Vert y_n-e_ny_ne_n\Vert_{E(\m,\tau)}=0.$$
In particular, $\{y_n\}_{n=1}^\infty$ has a subsequence that is
 equivalent to a disjointly supported basic 
sequence in $E$. Moreover, if $X$ has a basis then the sequence 
$\{y_n\}_{n=1}^\infty$ can be chosen to be a block basis of the basis of $X$.
\end{proposition}

\begin{proof} Assume that $X$ is not strongly embedded into $E(\m,\tau)$
and set $j: E(\m,\tau) \rightarrow \nmt$ the natural inclusion.
Since $X$ is not strongly embedded into $E(\m,\tau)$, the restriction
$j\vert_{X}$ is not an isomorphism. There exists a sequence 
$\{y_n\}_{n=1}^\infty$ in
the unit sphere of $X$ which converges to zero in measure.
Note that the bounded set $\{y_n, n\geq 1\}$ cannot be $E$-uniformly integrable.
Applying the non-commutative Kadec-Pe\l czy\'nski subsequence 
decomposition to the sequence $\{y_n\}_{n=1}^\infty$ on  $E(\m,\tau)$ (see 
\cite{Ran10} Theorem~3.7), there exist a 
subsequence of $\{y_n\}_{n=1}^\infty$ 
(which we will denote again by $\{y_n\}_{n=1}^\infty$ for simplicity) and a mutually
disjoint sequence of projections $\{e_n\}_{n=1}^\infty$ in $\m$ such that
the set
$\{y_n -e_ny_ne_n, n\geq 1\}$ is $E$-uniformly integrable and since 
$\{y_n -e_ny_ne_n\}_{n=1}^\infty$ converges to zero in measure, we get that
$\lim_{n\to \infty} \Vert y_n-e_ny_ne_n\Vert_{E(\m,\tau)}=0$.

Assume now that $X$ has a basis $\{x_n\}_{n=1}^\infty$. We will show that the
sequence $\{y_n\}_{n=1}^\infty$ above can be chosen to be a block basis
 of $\{x_n\}_{n=1}^\infty$.
In fact since $j(B_X)$ cannot be a neighborhood of zero for the (relative)
measure topology on $X$, for every $\epsilon>0$,
$B_{\mt}(0,\epsilon) \cap X \not\subset B_X$ (where $B_{\mt}(0,\epsilon)$ denotes
the ball centered at zero and with radius $\epsilon$ relative to the metric
of the measure topology). Denote by $\pi_n$ the projection $X$ onto 
$\overline{\text{span}}\{ x_k, k\leq n\}$.
Fix $z_1 \in S_X \cap B_{\mt}(0, 2^{-1})$ and choose $k_1\geq 1$ so that 
$\Vert z_1 - \pi_{k_1}(z_1)\Vert <2^{-1}$.
The restriction of $j$ on $(Id-\pi_{k_1})(X)$ cannot be an isomorphism.
As above, one can choose $z_2 \in S_X \cap B_{\mt}(0, 2^{-2})$ and 
$\pi_{k_1}(z_2)=0$. Inductively, one can construct a sequence 
$\{z_n\}_{n=1}^\infty$
in $S_X$ and a strictly increasing sequence of integers 
$\{k_n\}_{n=1}^\infty$ such that:
\begin{itemize}
\item[(i)] $z_n \in B_{\mt}(0, 2^{-n})$ for all $n \geq 1$;
\item[(ii)] $\Vert z_n -\pi_{k_n}(z_n) \Vert < 2^{-n}$ for all $n\geq 1$;
\item[(iii)]  $(Id-\pi_{k_n})(z_{n+1})=0$ for all $n\geq 1$.
\end{itemize}
Set $y_n:= \pi_{k_n}(z_n)$ for all $n\geq 1$. Clearly 
$\{y_n\}_{n=1}^\infty$ is a block basic sequence, $\Vert y_n \Vert^\alpha
 \geq 1-2^{-n\alpha}$ for
all $n\geq 1$ and  $\{y_n\}_{n=1}^\infty$ converges to zero in measure. 
The proof is complete.
\end{proof}

\smallskip

The next result can be viewed as a non-commutative analogue of
Proposition~1.c.10 of \cite{LT}(p.39). Below, $\{r_n(\cdot)\}_{n=1}^\infty$ 
denote the sequence of the Rademacher functions.

\begin{proposition}\label{dichotomy2}
Let $E$ be a symmetric Banach 
function space on $[0,\tau({\bf 1}))$. Assume that 
$E$ is order continuous and satisfy  the Fatou
property. Let $\{x_n\}_{n=1}^\infty$ be a sequence  
in $E(\cal M, \tau)$ with:
\begin{itemize}
\item[(i)] $\Vert x_n \Vert =1$ for all $n\geq 1$;
\item[(ii)] there exists a projection $e \in \cal{M}$ with
$\tau(e)<\infty$ and $ex_n=x_n$ for all $n\geq 1$.
\end{itemize}
Then either there exists a constant $C>0$ such that for every choice of
scalars $\{a_n\}_{n=1}^\infty$, we have
$$\int_{0}^1 \left\Vert \sum_{i=1}^n r_i(t)a_ix_i \right\Vert_{E(\m, \tau)}\ dt
\geq C \left(\sum_{i=1}^n \vert a_i \vert^2 \right)^{\frac{1}{2}}$$
for every $n\geq 1$, or
$\{x_n\}_{n=1}^\infty$ has a subsequence $\{x_{n_j}\}_{j=1}^\infty$ 
which is an unconditional basic
sequence equivalent to a disjoint element of $E$.
\end{proposition}

\begin{proof}
For $x \in E(\m,\tau)$, we set (as in \cite{SU}),
$$\sigma(x,\epsilon):= 
\ch_{[\epsilon\Vert x \Vert_{E(\m,\tau)}, \infty)}(\vert x \vert) $$
and
$$ M_{E(\m,\tau)}(\epsilon):=\{ x \in E(\m,\tau),\ 
\tau(\sigma(x,\epsilon))\geq \epsilon \}.$$
Assume first that for every $\epsilon >0$, there exists $n_\epsilon$ such that 
$\vert x_{n}^*\vert \notin M_{E(\m,\tau)}(\epsilon)$. We remark that
$\vert x_{n}^*\vert$ is supported by the finite projection $e$. There exists
a subsequence $\{x_{n_j}\}_{j=1}^\infty$ such that 
$\{\vert x_{n_j}^*\vert\}_{j=1}^\infty$ converges
to zero in measure. In particular, $\{x_{n_j}\}_{j=1}^\infty$ converge 
to zero in measure.
By the Kadec-Pe\l czy\'nski subsequence decomposition on $E(\m,\tau)$
(see \cite{Ran10}, Corollary~3.8),
 there exists a 
further subsequence (which we will denote again by $\{x_{n_j}\}_{j=1}^\infty$) and a
disjoint sequence of projections $\{e_j\}_{j=1}^\infty$ so that the set 
$\{x_{n_j} -e_j x_{n_j}e_j, j\geq 1\}$ is $E$-uniformly integrable so by
\cite{Ran10}(Proposition~2.8),
$$\lim_{j \to \infty}\Vert x_{n_j} -e_j x_{n_j}e_j \Vert =0$$
which shows that a subsequence of $\{x_{n_j}\}_{j=1}^\infty$ can be taken to be 
equivalent to a disjoint sequence of $E$.

On the other hand, if $\{\vert x_{n}^*\vert, n\geq 1\} \subset 
M_{E(\m,\tau)}(\epsilon)$ for some $\epsilon >0$ then
\begin{equation*}
\begin{split}
1=\Vert x \Vert &= \Vert \vert x_{n}^*\vert \Vert \cr
&\geq \Vert \vert x_{n}^*\vert \Vert_{L_1(\m,\tau) + \m} \cr
&\geq \left( \text{max}(1,\tau(e))\right)^{-1}
 \Vert \vert x_{n}^*\vert \Vert_{L_1(\m,\tau)} \cr
&\geq \epsilon \left( \text{max}(1,\tau(e))\right)^{-1}
 \tau\left(\sigma(\vert x_{n}^*\vert,\epsilon)\right) \cr
&\geq \epsilon^2 \left( \text{max}(1,\tau(e))\right)^{-1}
\end{split}
\end{equation*}
So for every $n\geq 1$, $\Vert x_n \Vert_1 = \Vert x_{n}^* \Vert_1 =
\Vert \vert x_{n}^*\vert \Vert_1 \geq 
\epsilon^2 \left( \text{max}(1,\tau(e))\right)^{-1}$.
Since $L_1(\m,\tau)$ is of cotype 2 (\cite{TJ2}), there exists $A_1>0$ such that
\begin{equation*}
\begin{split}
\int_{0}^1 \left\Vert \sum_{i=1}^n r_i(t)a_i x_i \right\Vert_{E(\m,\tau)}
 \ dt 
&=\int_{0}^1 \left\Vert e\left(\sum_{i=1}^n r_i(t)a_i x_i \right)
\right\Vert_{E(\m,\tau)} \ dt \cr
&\geq \int_{0}^1 \left\Vert e\left(\sum_{i=1}^n r_i(t)a_i x_i \right)
\right\Vert_{L_1(\m,\tau) + \m} \ dt \cr
&\geq \left( \text{max}(1,\tau(e))\right)^{-1}
\int_{0}^1 \left\Vert e\left(\sum_{i=1}^n r_i(t)a_i x_i \right)
\right\Vert_{L_1(\m,\tau)} \ dt \cr
&\geq A_1 \left( \text{max}(1,\tau(e))\right)^{-1}
\left(\sum_{i=1}^n \vert a_i \vert^2 \Vert x_i \Vert_{1}^2 \right)^{1/2} \cr
&\geq A_1 \epsilon^2 \left( \text{max}(1,\tau(e))\right)^{-2}
\left(\sum_{i=1}^n \vert a_i \vert^2 \right)^{1/2}.
\end{split}
\end{equation*}
The proof is complete.
\end{proof}

\begin{remarks}
We do not know if condition $(ii)$ can be removed. The same conclusion
holds if $(ii)$ is replaced by: there exists a projection
$e \in \m$ with $\tau(e)<\infty$ and  $x_ne=x_n$ for all $n\geq 1$.
\end{remarks}

\smallskip

The following theorem gives a description of symmetric basic sequences
in some symmetric space of measurable operators with type 2 and  can be 
viewed as semi-finite generalization of Theorem~2.4 of \cite{SU}.

\begin{theorem}\label{type2}
Let $E$ be an order continuous rearrangement invariant Banach 
function space on $[0,\tau({\bf 1}))$ with the Fatou
property and assume that  
 $E(\cal M, \tau)$ is of type 2. Then 
every symmetric basic sequence in $E(\m,\tau)$ either has a block
basic sequence equivalent to a disjointly supported sequence in $E$ or is
equivalent to $\ell_2$.
\end{theorem}

If $\tau({\bf 1})<\infty$, the theorem is a simple corollary of the above 
proposition with the word \lq\lq{block basic sequence}" replaced by 
\lq\lq{subsequence}".

For the general case,
choose a mutually orthogonal family $\{f_i\}_{i\in I}$
of projections in $\m$
with $\sum_{i\in I} f_i={\bf 1}$ for the strong operator topology and
$\tau(f_i)< \infty$ for all $i\in I$. Let $\{x_n\}_{n=1}^\infty$ be a symmetric 
basic sequence in $E(\m,\tau)$. Using a similar argument as in \cite{X}, one
can get an at most countable subset $\{f_k\}_{k=1}^\infty$ of 
$\{f_i\}_{i\in I}$
such that for each $f_i$ outside of $\{f_k\}_{k=1}^\infty$ and $n\geq 1$,
$f_ix_n=x_nf_i=0$. Let $f=\sum_{k \in\mathbb N} f_k$. For every $n\geq 1$,
we have
$fx_n=x_nf=x-n$. Replacing $\m$ by $f\m f$ and $\tau$ by its restriction on
$f\m f$, we may assume that $f={\bf 1}$. For every $n\geq 1$, set
$e_n:=\sum_{k=1}^n f_k$. The sequence $\{e_n\}_{n=1}^\infty$ is such that $e_n \uparrow_n {\bf 1}$
and $\tau(e_n) <\infty$ for all $n\geq 1$. 
Let $X:=\overline{\text{span}}\{x_n, n\geq 1\}$ and for $a\in \m$, let
$aX:=\{ax, x\in X\}$ and $Xa:=\{xa, x\in X\}$.

\begin{lemma}
If for every $n\geq 1$, $X$ is not isomorphic to $e_nX$ then 
there exists a normalized block basic sequence $\{y_k\}_{k=1}^\infty$ of 
$\{x_n\}_{n=1}^\infty$ and
a strictly  increasing sequence of integers $\{n_k\}_{k=1}^\infty$ so that
$$\Vert y_k -(e_{n_k}-e_{n_{k-1}})y_k \Vert<2^{-k},\ \ \ \text{for}\ 
k\geq 1.$$
Similarly, if for every $n\geq 1$, $X$ is not isomorphic to $Xe_n$ then 
there exists a normalized block basic sequence $\{y_k\}_{k=1}^\infty$ of 
$\{x_n\}_{n=1}^\infty$ and
a strictly increasing sequence of integers $\{n_k\}_{k=1}^\infty$ so that
$$ \Vert y_k -y_k(e_{n_k}-e_{n_{k-1}}) \Vert <2^{-k},\ \ \ \text{for}\ 
k\geq 1.$$ 
\end{lemma}
\begin{proof}
Inductively, we will construct a sequence $\{y_k\}_{k=1}^\infty$ in the unit sphere
of $X$, strictly increasing sequences of integers 
$\{m_k\}_{k=1}^\infty$ and $\{n_k\}_{k=1}^\infty$  such
that:
\begin{itemize}
\item[(i)] $y_k \in \text{span}\{x_n, m_{k-1}<n \leq m_k\}$ for all 
$k \geq 1$;
\item[(ii)] $\Vert e_{n_{k_1}}y_k \Vert < 2^{-(k+1)}$ for all $k \geq 1$;
\item[(iii)] $\Vert y_k - e_{n_{k}}y_k \Vert < 2^{-(k+1)}$ for all $k \geq 1$.
\end{itemize}
Fix $y_1$ a finitely supported vector in $S_X$ and let $m_1 \geq 1$ so that
$y_1 \in \text{span}\{x_n, n \leq m_1\}$. Since 
$({\bf 1}- e_{n}) \downarrow_n 0$, there exists $n_1$ such that 
$\Vert y_1 -e_{n_1}y_1 \Vert <2^{-1}$.

Assume that the construction is done for $1, 2, \cdots, (j-1)$. Let
$X_j=\overline{\text{span}}\{x_n, n \geq m_{j-1}\}$. Since $X_j$ is not
isomorphic to $e_{n_{j-1}}X_j$, there exist $y_j \in S_{X_j}$ such that
$\Vert e_{n_{j-1}}y_j\Vert < 2^{-(j+1)}$. 
By perturbation, we can assume that $y_j$ is finitely supported. If we fix
$n_j> n_{j-1}$ so that $\Vert y_j - e_{n_j}y_j\Vert < 2^{-(j+1)}$ then
$\Vert y_j - (e_{n_j}- e_{n_{j-1}})y_j\Vert < 2^{-j}$ and the lemma follows.
\end{proof}

To prove the theorem, assumme first the there exist $n_0 \geq 1$ such that
$X$ is isomorphic to $e_{n_0}X$. Since $\tau(e_{n_0})<\infty$, the sequence 
$\{e_{n_0}x_n\}_{n=1}^\infty$ satisfies the assumptions of 
Proposition~\ref{dichotomy2}. Since
 $E(\m,\tau)$ has type 2, either $\{e_{n_0}x_n\}_{n=1}^\infty$ is equivalent to $\ell_2$
or there exists a subsequence $\{e_{n_0}x_{n_j}\}_{j=1}^\infty$  which is equivalent 
to a sequence of disjoint elements of $E$ and by isomorphism, the theorem follows.

Assume now that for every $n \geq 1$, $X$ is not isomorphic to 
$e_nX$. By the above lemma, there exist a normalized 
block basic sequence $\{y_k\}_{k=1}^\infty$
and a strictly increasing sequence of integers $\{n_k\}_{k=1}^\infty$
so that for every $k \geq 1$,
\begin{equation}\label{2.1}
\Vert y_k - (e_{n_k}- e_{n_{k-1}})y_k\Vert < 2^{-k}.
\end{equation}
Let $Y:=\overline{\text{span}}\{(e_{n_k}- e_{n_{k-1}})y_k, k\geq 1\}$.
As above, if there exists $m_0$ such that $Y$ is isomorphic to $Ye_{m_0}$,
then the conclusion follows. Otherwise,
there exist a block basic sequence $\{z_k\}_{k=1}^\infty$ of 
$\{(e_{n_k}- e_{n_{k-1}})y_k\}_{k=1}^\infty$ and a
strictly increasing sequence of 
integers $\{m_k\}_{k=1}^\infty$ such that for every $k\geq 1$,
\begin{equation}\label{2.2} 
\Vert z_k -z_k(e_{m_k}- e_{m_{k-1}})\Vert <2^{-k}.
\end{equation}
We remark that since the sequence $\{z_k\}_{k=1}^\infty$ is a block basic sequence of
$\{(e_{n_k}- e_{n_{k-1}})y_k\}_{k=2}^\infty$, there exists a sequence 
$\{q_k\}_{k=1}^\infty$ of
mutually disjoint projections  such that for every $k\geq 1$,
$z_k=q_kz_k$. Therefore, the sequence
$\{z_k(e_{m_k}- e_{m_{k-1}})\}_{k=2}^\infty$  is both right and left
disjointly supported and hence is equivalent to a disjointly supported
sequence in $E$. By (\ref{2.2}), we conclude that $\{z_k\}_{k=1}^\infty$ 
has a subsequence that is equivalent 
to a disjointly supported sequence in $E$ (see for instance,
\cite{D1}, Theorem~9 p.46). Since 
$\{z_k\}_{k=1}^\infty$ is a block basic
sequence of $\{(e_{n_k}- e_{n_{k-1}})y_k\}_{k=2}^\infty$, inequality (\ref{2.1}) 
shows
that the corresponding  block  of $\{y_k\}_{k=1}^\infty$  is equivalent to
$\{z_k\}_{k=1}^\infty$. The proof of the theorem is complete.
\qed

\section{Subspaces of Lorentz spaces}

In this section, we will specialize to the concrete case of Lorentz
spaces. We begin by recalling some definitions and 
basic facts about Lorentz spaces.

For $0<p<\infty$, $0<q<\infty$, and $I = [0,1]$ or $[0,\infty)$, the Lorentz
function space $L_{p,q}(I)$ is the space of all (classes of) Lebesgue
measurable functions $f$ on $I$ for which $\Vert f \Vert_{p,q} < \infty$,
 where
\begin{equation}\label{definition}
\begin{split}
\Vert f \Vert_{p,q} &= \left(\int_{I} \mu_{t}^{q}(f)\ 
d(t^{q/p})\right)^{1/q}, \quad q < \infty, \cr
&= \sup_{t \in I} t^{1/p} \mu_t(f), \quad q= \infty.
\end{split}\end{equation}
Clearly, $L_{p,p}(I)=L_p(I)$ for any $p>0$. It is well known that
for $1 \leq q \leq p < \infty$, (\ref{definition}) defines a norm under which 
$L_{p,q}(I)$ is a separable rearrangement invariant Banach function space;
otherwise, (\ref{definition}) defines a quasi-norm on $L_{p,q}(I)$ ( which is known to be 
equivalent to a norm  
 if $1<p<q<\infty$). The following lemma was observed in \cite{CADI2}. It
contains the technical ingredients for the construction of the 
non-commutative counterparts. 

\begin{lemma}
Let $0<p<\infty$, $0<q<\infty$.
\begin{itemize}
\item[(i)] If $q<p$, then $L_{p,q}(I)$ is $q$-convex with constant $1$ and
satisfies a lower $p$-estimate with constant $1$;
\item[(ii)] $L_{p,q}(I)$ satisfies an upper $r$-estimate
and  lower $s$-estimate (with some constant $C$) where 
$r=\text{min}(p,q)$ and $s=\text{max}(p,q)$.
\end{itemize}
\end{lemma}

For $0<p<q<\infty$, $L_{p,q}(I)$ can be equivalently renormed to be a 
quasi-Banach lattice that is $\gamma$-convex (for $\gamma <p$) with 
constant $1$ and satisfies a lower $q$-estimate of constant $1$. Hence
for any $0<p,q <\infty$, we can define the non-commutative 
space $L_{p,q}(\m,\tau)$ as in Section~2. 
Since we are only interested in isomorphic properties, we will
use the quasi-norm defined in (\ref{definition}). All results
from sectio~2 apply to $L_{p,q}(\m,\tau)$ with appropriate
values of $p$ and $q$.

  For any given $0<p<\infty$ and $0<q\leq \infty$,
it is well known that the space
 $L_{p,q}(I)$  is equal
(up to an equivalent quasi-norm) to the spaces $(L_{p_{1}}(I),
L_{p_{2}}(I))_{\theta,q}$ constructed using the real interpolation method
where $0<p_{1}<p_{2}<\infty$,  $0 < \theta < 1$ and $1/p=
(1-\theta)/{p_{1}}+ \theta/{p_{2}}$. This was extended to the
non-commutative setting by Xu (\cite{X2}).
\begin{lemma}\label{interpolation}
For $0<p_1, p_2, q<\infty$ and 
 $0<\theta<1$ then
$$(L_{p_1}(\m,\tau), L_{p_2}(\m,\tau))_{\theta,q} =
L_{p,q}(\m,\tau)$$
(with equivalent quasi-norms) where 
$1/p=(1-\theta)/{p_{1}}+ \theta/{p_{2}}$.
\end{lemma}

The following result can be found in \cite{CADI}(Lemma~2.4 and the 
remark after Theorem~2.5).

\begin{lemma}
Let $0<p<\infty$ and $0< q<\infty$. Let $\{f_n\}_{n=1}^\infty$ be a 
normalized
disjointly supported sequence in $L_{p,q}[0,\infty)$.
Then $\overline{\text{span}}\{f_n, n\geq 1\}$ contains a copy
of $\ell_q$.
\end{lemma}

A non-commutative extension follows directely from 
Proposition~\ref{disjoint}.

\begin{proposition}\label{lq}
Let $0<p<\infty$ and $0<q<\infty$. 
Let $\{x_n\}_{n=1}^\infty$ be a normalized basic 
 sequence in $L_{p,q}(\m,\tau)$. If $\{x_n\}_{n=1}^\infty$ is both 
right and left disjointly supported then 
 $\overline{\text{span}}\{x_n, n\geq 1\}$ contains a copy
of $\ell_q$.
\end{proposition}

The next result is an analogue of Theorem~2.5 of \cite{CADI}
and is a direct application of Proposition~\ref{dichotomy1} and 
Proposition~\ref{lq}

\begin{theorem}\label{dichotomy3}
Let  $0<p<\infty$, $0< q<\infty$, and let $X$ be a subspace
of $L_{p,q}(\m,\tau)$. Then either $X$ is strongly embedded 
into $L_{p,q}(\m,\tau)$  or $X$ contains a copy of $\ell_q$.
\end{theorem}

\smallskip

For the next result, we need to fix some notation.
Let $\cal N$ be a von Neumann algebra on a given Hilbert space $H$ with
semi-finite trace $\varphi$.  Define $[\cal N]$ to be the von Neumann
algebra over $\ell_2(H)$ as follows:
$$
[\cal N]: = \left\{(a_{ij})_{ij};\  \forall \ i,j,\  a_{ij} \in \cal N,  \Vert
(a_{ij})_{ij}\Vert_{B(\ell^{2}(H))} < \infty \right\}.
$$
equipped with the trace 
$[\varphi]((a_{ij})_{ij}) = \sum_{i=1}^\infty \varphi(a_{ii})$. It is clear
that 
$([\cal N],[\varphi])$ is a semi-finite von Neumann algebra acting on 
$\ell_2(H)$.
Let $\{y_{k}\}_{k=1}^\infty$  be a sequence in $\cal N$.  
For each $k\geq 1$, we define $[y_{k}]=([y_{k}]_{ij})_{ij}$ 
by setting:
$[y_{k}]_{1,k} = y_{k}$ and $[y_{k}]_{ij} =0$ for $(i,j) \neq (1,k)$
i.e for $k\geq 1$,
$$
[y_{k}] :=
\pmatrix
0 & \hdots & 0 & y_{k} & 0 & \hdots  \\
0 & \hdots & 0 & 0    & 0 & \hdots  \\
\vdots & \hdots  & \vdots & \vdots & \vdots &  \hdots  
\endpmatrix.
$$
This amounts to placing the sequence $\{y_k\}_{k=1}^\infty$ in the first
 row of an infinite matrix.
\begin{lemma}\label{khintchine}
Let $0<p<2$ and 
 $\{y_{k}\}_{k=1}^\infty$ be a sequence in $L_{p,q}(\cal N, \varphi)$. 
 There
exists an absolute constant $C$ such that for every choice  of
scalars $\{a_k\}_{k=1}^\infty$,
\begin{equation}
%\begin{split}
\int^{1}_{0} \left\Vert \sum^{n}_{k=1} r_{k}(t) a_{k} y_{k} 
 \right\Vert_{L_{p,q}(\cal N, \varphi)}^2 \ dt 
\leq
C \text{min}\left\{ \left\Vert \sum^{n}_{k=1}  a_{k} 
[y_{k}]
\right\Vert_{L_{p,q}([\cal N], [\varphi])}^2,
 \left\Vert \sum^{n}_{k=1} \bar{a}_{k}[y^{*}_{k}]
\right\Vert_{L_{p,q}([\cal N],[\varphi])}^2\right\}
%\end{split}
\end{equation}
for every $n\geq 1$.
\end{lemma}

\begin{proof}
We claim first that for $0<p<2$, there exists an absolute constant $A$ 
such that:
\begin{equation}\label{Kin}
\left(\int_{0}^1 \left\Vert \sum_{k=1}^n  r_{k} (t) a_{k} y_{k} 
\right\Vert_{L_{p} (\cal N,\varphi)}^2 \ dt\right)^{\frac{1}{2}}  \leq 
\sqrt{A}
\left\Vert (\sum^{n}_{k=1} \vert a_{k} \vert^{2} y_k y^{*}_{k} 
)^{\frac{1}{2}} \right\Vert_{L_{p}(\cal N,\varphi)}.
\end{equation}
Indeed,
for $1 \leq p < 2$, (\ref{Kin}) follows from non-commutative Khintchine inequality
(see\cite{LPI} Corollary III-4 and Remark III-6).  For $0 < p < 1$, we
remark that for any given finite sequence $\{x_{k}\}_{k=1}^ n$ in
$ L_{p}(\cal N,\varphi)$, 
$$
\sum^{n}_{k=1} \Vert x_{k} \Vert^{2}_{L_{p}(\cal N,\varphi)}  
\leq  A \Vert
(\sum_{k=1}^n x_{k} x_{k}^{*})^{\frac{1}{2}} \Vert_{L_{p}(\cal N,\varphi)}
$$
for some fixed constant $A$.  We observe that
\begin{equation*}
\begin{split}
\vert \sum_{k=1}^n r_{k}(t) [x_{k}]^{*} \vert^{2} &=
\left(\sum_{k=1}^n r_k(t)[y_{k}]\right). \left(\sum_{k=1}^n 
 r_{k}(t)[x_{k}]^{*}\right) \cr
&=
\matrix
\pmatrix
r_{1}(t) x_{1} & r_{2}(t)x_{2} & \hdots & r_{n}(t)x_{n} & \hdots  \\
0           &   0        & \hdots & 0          & \hdots  \\
\vdots      &   \vdots   &  \hdots      &  \vdots &\hdots
\endpmatrix
\cdot
& \pmatrix
r_{1}(t) x^{*}_{1} & 0 & \hdots  \\
r_{2}(t) x^{*}_{2}       & 0 & \hdots  \\
\vdots                & \vdots &\hdots
\endpmatrix
\endmatrix
\cr
&=
\pmatrix
\sum_{k=1}^n   x_{k} x^{*}_{k} & 0 & \hdots \\
0  &0  &\hdots\\
\vdots  &\vdots &\hdots 
\endpmatrix 
\end{split}
\end{equation*}
so
$$
\vert \sum_{k=1}^n  r_{k}(t) [x_{k}]^{*} \vert^{p} =
\pmatrix
\left(\sum_{k=1}^n  x_{k} x^{*}_{k}
 \right)^{\frac {p}{2}} & 0
& \hdots  \\
0 & 0 & \hdots &  \\
\vdots &   \vdots &  \hdots  \
\endpmatrix
$$
which implies that
$$
\left\Vert (\sum_{k=1}^n x^{*}_{k}x_{k})^{\frac{1}{2}} 
\right\Vert^{2}_{L_{p}(\cal N,\varphi)}
= \int_{0}^1 \left\Vert \sum^{n}_{k=1}
r_{k}(t)[x_{k}]^{*} \right\Vert^{2}_{L_{p}([\cal N],[\varphi])}\ dt
$$
and since $L_{p}([\cal N],[\varphi])$ is of cotype $2$(\cite{X2}),
there exists a constant $A>0$ so that,
\begin{equation*}
\begin{split}
A \left\Vert (\sum_{k=1}^n x_k x^{*}_{k})^{\frac{1}{2}} 
\right\Vert^{2}_{L_{p}(\cal N,\varphi)}
&\geq \sum^{n}_{i=1} \left\Vert [x_{k}]^{*} 
\right\Vert^{2}_{L_{p}([\cal N],[\varphi])}
\cr
&= \sum^{n}_{k=1} \left\Vert x_{k} \right\Vert^{2}_{L_{p}(\cal N,\varphi)}.
\end{split}
\end{equation*}
For each $\theta \in \{-1,1\}^n$, set
$ s_\theta=\sum \limits^{n}_{k=1} \theta_{k} a_ky_{k}$. Applying the 
above inequality to the finite sequence $\{s_\theta\}_\theta$, we get   
 that
\begin{equation*}
\begin{split}
\int^{1}_{0} \left\Vert \sum^{n}_{k=1} r_{k} (t) a_ky_{k} 
\right\Vert_{L_p(\cal N, \varphi)}^{2} dt &= 2^{-n}
\sum \limits_{\theta \in \{-1,1\}^n} \left\Vert s_{\theta} 
\right\Vert_{L_p(\cal N, \varphi)}^{2} \cr
&\leq 2^{-n} A \left\Vert (\sum \limits_{\theta \in \{-1,1\}^n} 
s_{\theta} s^{*}_{\theta} )^{\frac{1}{2}} \right\Vert_{L_p(\cal N, \varphi)}^{2} \cr
&= A 2^{-n} \left\Vert (\sum \limits_{\theta \in \{-1,1\}^n} \sum\limits_{i,k} 
 \theta_{i} 
\theta_{k} \bar{a_i} a_k y_ky^{*}_{i} )^{\frac {1}{2}} 
\right\Vert_{L_p(\cal N, \varphi)}^{2} \cr
&= A \left\Vert (\int_{0}^1 \sum\limits_{i,k} r_{k}(t) r_{i}(t) 
\bar{a_i}a_k y_ky_{i}^{*} )^{\frac{1}{2}} \right\Vert_{L_p(\cal N, \varphi)}
\cr
&= A \left\Vert (\sum_{k=1}^n  
\vert a_k  \vert^2 y_k y_{k}^{*} )^{\frac{1}{2}} 
\right\Vert_{L_p(\cal N, \varphi)}^{2}
\end{split}
\end{equation*}
and therefore (\ref{Kin}) is verified.

For the proof of the lemma, we note as above that
$$
\vert \sum_{k=1}^n  \bar{a}_{k} [y_{k}]^{*} \vert^{p} =
\pmatrix
\left(\sum_{k=1}^n \vert a_{k} \vert^{2} y_{k} y^{*}_{k}
 \right)^{\frac {p}{2}} & 0
& \hdots  \\
0 & 0 & \hdots &  \\
\vdots &   \vdots &  \hdots  \
\endpmatrix
$$
therefore
$$
 \left\Vert \sum_{k=1}^n  \bar{a}_{k} [y_{k}]^{*} 
\right\Vert_{L_{p}([\cal N], [\varphi])} =
\left\Vert (\sum_{k=1}^n \vert a_{k} \vert^{2} y_{k} y_{k}^{*})^{\frac {1}{2}}
\right\Vert_{L_{p}(\cal N, \varphi)},
$$
hence
$$
 \left\Vert \sum_{k=1}^n  {a}_{k} [y_{k}] 
\right\Vert_{L_{p}([\cal N], [\varphi])} =
\left\Vert (\sum_{k=1}^n \vert a_{k} \vert^{2} y_{k} y_{k}^{*})^{\frac {1}{2}}
\right\Vert_{L_{p}(\cal N, \varphi)}
$$
and by (\ref{Kin}),
$$
\int^{1}_{0} \left\Vert \sum^{n}_{k=1} r_{k}(t) a_{k} y_{k} 
 \right\Vert_{L_{p}(\cal N, \varphi)}^2 \ dt 
\leq
A \left\Vert \sum^{n}_{k=1}  a_{k} 
[y_{k}]
\right\Vert_{L_{p}([\cal N], [\varphi])}^2.
$$

\begin{sublemma}
 For every $0 < p < 1$, the
map $(a_{ij})_{ij} \longrightarrow \sum \limits_{k} r_k(\cdot) a_{1k}$
is bounded as a linear map from $L_{p}([\cal N],[\varphi])$ into
$L_{2}([0,1], L_{p}(\cal N,\varphi))$.
\end{sublemma}
Let $a=(a_{ij})_{ij}$ be an element of $L_{p}([\cal N], [\varphi])$ and 
consider $\vert a^*\vert^2 =(b_{ij})_{ij}$. Clearly,
$b_{11}= \sum_{k=1}^\infty a_{1k}a_{1k}^*$. Let $e$ be the projection in
$[\cal N]$ defined by $e=(e_{ij})_{ij}$ with $e_{11}= {\bf 1}$ and
$e_{ij}=0$ for $(i,j)\neq (1,1)$. We have
$$
e\vert a^*\vert^2 e=
\pmatrix
\sum_{k=1}^\infty  a_{1k} a^{*}_{1k}
  & 0
& \hdots  \\
0 & 0 & \hdots &  \\
\vdots &   \vdots &  \hdots  
\endpmatrix
$$
so $\Vert e\vert a^*\vert^2 e\Vert_{L_{p/2}([\cal N], [\varphi])}
=\Vert (\sum_{k=1}^\infty a_{1k}a_{1k}^*)^{1/2}
\Vert_{L_{p}(\cal N, \varphi)}^2$
and as above,
\begin{equation*}
\begin{split}
\int^{1}_{0} \left\Vert \sum^{\infty}_{k=1} r_{k}(t) a_{1k} 
 \right\Vert_{L_{p}(\cal N, \varphi)}^2 \ dt 
&\leq
A \left\Vert \sum^{\infty}_{k=1}   
[a_{1k}]
\right\Vert_{L_{p}([\cal N], [\varphi])}^2 \cr
&=A\left\Vert (\sum_{k=1}^\infty a_{1k}a_{1k}^*)^{1/2}
\right\Vert_{L_{p}(\cal N, \varphi)}^2 \cr
&=A\left\Vert e\vert a^* \vert^2 e 
\right\Vert_{L_{p/2}([\cal N], [\varphi])} \cr
&\leq A\left\Vert a \right\Vert_{L_{p}([\cal N], [\varphi])}^2.
\end{split}
\end{equation*}
and the sublemma follows.

By interpolation, it is also a bounded map from $L_{p,q}([\cal N], [\varphi])$ into
$L_{2}([0,1], L_{p,q}(\cal N, \varphi))$ which shows in particular that there exists an
absolute constant $C$ such that
$$
\int_{0}^1 \left\Vert \sum_{k=1}^n r_k(t) a_{k} y_{k} 
\right\Vert_{L_{p,q}(\cal N,\varphi)}^2\ 
dt \leq C \left\Vert \sum_{k=1}^n  a_{k} [y_{k}] 
\right\Vert_{L_{p,q}([\cal N],[\varphi])}^2.
$$
By taking adjoints, the other inequality follows. The proof of 
the lemma is complete.
\end{proof}

Our next result is a disjointification of sequences in $L_{p,q}(\cal M,\tau)$
and could be of independent interest.

\begin{proposition}\label{disjointification}
Let $(\m,\tau)$ be a semi-finite von Neumann algebra. There exists a 
semi-finite von Neumann algebra $\cal S$ equipped with a faithful
normal semi-finite trace $\omega$ with the following properties:
\begin{itemize}
\item[(i)] $\m $ is a von Neumann subalgebra of $\cal S$;
\item[(ii)] $\tau$ is the restriction of $\omega$ on $\m$;
\item[(iii)] For $0<p< 2$ and $0<q<\infty$, there exists a
constant $K$ such that for any given basic sequence
$\{x_n\}_{n=1}^\infty$ in $L_{p,q}(\m,\tau)$, there exists a left
and right
disjointly supported  sequence $\{s_n\}_{n=1}^\infty$ in $L_{p,q}(\cal S, \omega)$
such that for any choice of scalars $\{a_k\}_{k=1}^\infty$ and
$n \geq 1$,
$$ 
\int_{0}^1 \left\Vert \sum_{k=1}^n a_k r_k(t)x_k 
\right\Vert_{L_{p,q}(\m,\tau)}^2\ dt \leq
K\left\Vert \sum_{k=1}^n a_ks_k 
\right\Vert_{L_{p,q}(\cal S,\omega)}^2.
$$
\end{itemize}
\end{proposition} 

\begin{proof}
Using the notation above, let $\cal{N}=[\m]$, $\varphi =[\tau]$.
Clearly, $(\cal N, \varphi)$ is a semi-finite von Neumann algebra on
the Hilbert space $H=\ell_2(\cal H)$. 
Set $\cal{S}=[\cal{N}]$ and $\omega=[\varphi]$. As above, 
$\m$ can be identified as a von Neumann subalgebra of $\cal S$ with 
$\tau$ being the restriction of $\omega$ on $\m$.

Let  $\{x_{n}\}_{n=1}^\infty$  be a basic sequence  in
$L_{p,q}(\m, \tau)$. Consider the sequence $\{[x_n]\}_{n=1}^\infty$
in $\cal{N}=[\m]$.

\noindent
{\it Claim:\ The sequence $\{[x_n]\}_{n=1}^\infty$ is right disjointly
supported.}

To verify this claim, recall that elements of $\cal N$ are infinite
matrices with entries in $\m$. For $n\geq 1$, let
$\pi_n =(a_{ij})_{ij}$ with $a_{n,n}={\bf 1}$ and $a_{i,j}=0$ for
$(i,j)\neq (n,n)$. Clearly $\{\pi_n\}_{n=1}^\infty$is a mutually 
disjoint sequence of projection in $\cal N$ and for each
$n\geq 1$, $[x_n]\pi_n = [x_n]$.

For each $n\geq 1$, let $z_n=[x_n] \in L_{p,q}(\cal{N}, \varphi)$ 
and consider the sequence
$\{s_n \}_{n=1}^\infty$ in $L_{p,q}(\cal S, \omega)$ defined
by $s_n:=[z_{n}^*]^*$.

\noindent
{\it Claim:\ The sequence $\{s_n\}_{n=1}^\infty$ is left and right
disjointly supported.}

First note  that, as  above, the sequence $\{[z_{n}^*]\}_{n=1}^\infty$ 
is right disjointly supported  so its adjoints $\{s_n\}_{n=1}^\infty$
is left disjointly supported. To prove that it is right disjointly supported,
consider the following sequence $\{e_n\}_{n=1}^\infty$ in $\cal S$: 
$e_n=(a_{ij}^{(n)})_{ij}$ where $a_{11}^{(n)}=\pi_n $ and 
$a_{ij}^{(n)}=0$ for $(i,j)\neq (1,1)$.

It is clear that the $e_n$'s are projections in $\cal S$ and since 
$\{\pi_n\}_{n=1}^\infty$ is mutually disjoint in $\cal{N}$,
$\{e_n\}_{n=1}^\infty$ is mutually disjoint and one can see that for every
$n\geq 1$, $s_ne_n=s_n$.

To complete the proof, we use  Lemma~\ref{khintchine},
\begin{equation*}\begin{split}
\int_{0}^1 \left\Vert \sum_{k=1}^n r_k(t) a_{k} x_{k} \right\Vert_{L_{p,q}(\m,\tau)}^2
\ dt
&\leq C\left\Vert \sum_{k=1}^n a_k [x_k] 
\right\Vert_{L_{p,q}(\cal N, \varphi)}^2 \cr
&=C\int_{0}^1 \left\Vert \sum_{k=1}^n r_k(t)a_k [x_k] 
\right\Vert_{L_{p,q}(\cal N, \varphi)}^2 \ dt \cr
&=C\int_{0}^1 \left\Vert \sum_{k=1}^n r_k(t)a_k z_k 
\right\Vert_{L_{p,q}(\cal N, \varphi)}^2 \ dt.
\end{split}
\end{equation*}
Applying Lemma~\ref{khintchine} on the von Neumann algebra $\cal N$,
\begin{equation*}
\begin{split}
\int_{0}^1 \left\Vert \sum_{k=1}^n r_k(t) a_{k} x_{k} \right\Vert_{L_{p,q}(\m,\tau)}^2
\ dt
&\leq C^2 \left\Vert \sum_{k=1}^n \bar{a_k} [z_{k}^*] 
\right\Vert_{L_{p,q}(\cal S, \omega)}^2  \cr
&=C^2 \left\Vert \sum_{k=1}^n {a_k} [z_{k}^*]^* 
\right\Vert_{L_{p,q}(\cal S, \omega)}^2  \cr
&=C^2 \left\Vert \sum_{k=1}^n {a_k} s_k 
\right\Vert_{L_{p,q}(\cal S, \omega)}^2.  
\end{split}\end{equation*}
The proof is complete
\end{proof}

\smallskip

We may now state the main result of this paper.

\begin{theorem}
Let $0<p<\infty$, $0<q<\infty$, $p\neq q$ and $p\neq 2$. Then
$\ell_p$ does not embed into $L_{p,q}(\m,\tau)$. In particular,
the Lorentz-Schatten  ideals $S_{p,q}$ does not contain $\ell_p$.
\end{theorem}

The proof will be divided into several cases. First, notice that since 
$p\neq q$, Proposition~\ref{dichotomy3} shows that every subspace of
$L_{p,q}(\m,\tau)$ equivalent to $\ell_p$ (and therefore
not containing any copy of $\ell_q$) is strongly embedded into $L_{p,q}(\m,\tau)$.
Fix $r>q$ then $\Vert \cdot \Vert_{p,r} \leq C\Vert \cdot \Vert_{p,q}$,
where $C$ is a constant depending only on $p$, $q$ and $r$ (see
for instance \cite{BENSHA}, Proposition~4.2 p.217). In particular, 
there exists a continuous inclusion from $L_{p,q}(\m,\tau)$ into
$L_{p,r}(\m,\tau)$ and if $X$ is a strongly embedded subspace of 
$L_{p,q}(\m,\tau)$ then $X$ is isomorphic to a subspace of 
$L_{p,r}(\m,\tau)$ so without loss of generality, we can assume that
$p<q$ and $1<q$.

\smallskip

\noindent
{\it Case $0<p<q<\infty$ and $p<2$.}

Assume that  there exists a sequence $\{x_n\}_{n=1}^\infty$ that
is $M$-equivalent to $\ell_p$ in $L_{p,q}(\m,\tau)$ and consider the 
disjoint sequence
$\{y_n\}_{n=1}^\infty$ in $L_{p,q}(\cal S,\omega)$ as in 
Proposition~\ref{disjointification}.
For every finite sequence of scalars $\{a_n\}$, we have:
\begin{equation*}
\begin{split}
\left(\sum_n \vert a_{n} \vert^{p}\right)^{\frac{1}{p}} 
&\leq M \left(\int_{0}^1 \left\Vert \sum_n
 r_n(t) a_{n} x_{n} \right\Vert_{L_{p,q}(\cal M, \tau)}^2  \ dt 
\right)^{\frac12} \cr
&\leq M.\sqrt{K}  \left\Vert\sum_n  a_{n} y_n
\right\Vert_{L_{p,q}(\cal S,\omega)}  \cr
&\leq N.M.\sqrt{K}\left\Vert \sum_n a_{n} \varphi_{n}
\right\Vert_{L_{p,q}(0,\infty)}  
\end{split}
\end{equation*}
where $\{\varphi_{k}\}_{k=1}^\infty$ is a disjoint sequence in 
$L_{p,q}[0,\infty)$ and $N>0$. Since $p<q$, the space $L_{p,q}[0,\infty)$
satisfies an upper $p$-estimate hence there exists constants $C_1$ and
$C_2$ such that
$$C_1\left( \sum_{k=1}^n \vert a_k \vert^p \right)^{\frac{1}{p}}
\leq \left\Vert \sum_{k=1}^n a_k \varphi_k 
\right\Vert_{L_{p,q}[0,\tau({\bf 1}))}
\leq C_2\left( \sum_{k=1}^n \vert a_k \vert^p \right)^{\frac{1}{p}}.$$
But this is a contradiction since 
$\overline{\text{span}}\{\varphi_k, k\geq 1 \}$ contains a copy of 
$\ell_q$.
\qed

\smallskip

\noindent
{\it Case  $2<p<q<\infty$.}

 We remark that $L_{p,q}(\m,\tau)$ is of type~2 (see for instance,
\cite{LT}, Proposition~2g.22 p.230).
Assume that there exists a sequence $\{x_n\}_{n=1}^\infty$ in
 $L_{p,q}(\m,\tau)$ that
is equivalent to $\ell_p$. Since $p\neq 2$, Theorem~\ref{type2} 
 and 
Proposition~\ref{disjoint} imply that 
$\{x_n\}_{n=1}^\infty$ contains a block basic sequence $\{y_n\}_{n=1}^\infty$ 
that is equivalent to a disjointly supported normalized sequence in
$L_{p,q}[0,\tau({\bf 1}))$ so 
$\overline{\text{span}}\{y_n, n\geq 1\}$ does not contain $\ell_p$. 
 Contradiction.

\qed

\smallskip

The following application characterizes strongly embedded subspaces in
$L_p(\m,\tau)$. It generalizes results of Rosenthal and Kalton on
$L_p[0,1]$.

\begin{corollary}\label{Rosenthal}
Let $0<p<\infty$, $p\neq 2$ and $X$ be a subspace of $L_p(\m,\tau)$. Then the 
following are equivalent:
\begin{itemize}
\item[(1)] $X$ contains $\ell_p$;
\item[(2)] $X$ is not strongly embedded into $L_p(\m,\tau)$.
\end{itemize}
\end{corollary}

\begin{proof}
Let $X$ be a subspace of $L_p(\m,\tau)$ and assume that $X$ contains
$\ell_p$. Since for $p<q$, as above,
$\Vert \cdot \Vert_{p,q} \leq C\Vert \cdot \Vert_p$, for some constant $C$
(see  
\cite{BENSHA} Proposition~4.2  p.217). There exists an inclusion map
from $L_p(\m,\tau)$ into $L_{p,q}(\m,\tau)$. If $X$ is strongly embedded
into $L_p(\m,\tau)$, then $X$ is isomorphic to a subspace of 
$L_{p,q}(\m,\tau)$.
In particular $\ell_p$ embeds into $L_{p,q}(\m,\tau)$. Contradiction.

The converse is a direct consequence of Theorem~\ref{dichotomy3}.
\end{proof}

The next result is known for copies of $\ell_1$ in preduals
of von Neumann algebras.

\begin{corollary}
Let $1\leq p <\infty$, $p\neq 2$. If $\{x_n\}_{n=1}^\infty$ is a sequence in $L_p(\m,\tau)$
that is equivalent to $\ell_p$ and $\{\varepsilon_n\}_{n=1}^\infty$ 
is a sequence in
the interval $(0, 1)$ with $\varepsilon_n \downarrow_n 0$, then there 
exists a block basis  $\{y_n\}_{n=1}^\infty$ of $\{x_n\}_{n=1}^\infty$ 
such that:
\begin{equation*}
\left(\sum_n \vert a_n \vert^p\right)^{\frac{1}{p}} -
\left(\sum_n \vert a_n \vert^p {\varepsilon_n}^p\right)^{\frac{1}{p}}
\leq \left\Vert \sum_n a_ny_n \right\Vert
\leq \left(\sum_n \vert a_n \vert^p\right)^{\frac{1}{p}} +
\left(\sum_n \vert a_n \vert^p {\varepsilon_n}^p\right)^{\frac{1}{p}}
\end{equation*}
for all finite sequence $(a_n)_n$ of scalars. In particular,
for every $k\geq 1$,  the sequence 
$\{y_n\}_{n=k}^\infty$ is $(1+\varepsilon_k)$-equivalent to $\ell_p$.
\end{corollary}
\begin{proof}
Since $\ell_p$ is not strongly embedded into $L_p(\m,\tau)$,
Proposition~\ref{dichotomy1} implies the existance of a block basic sequence
$\{z_n\}_{n=1}^\infty$ of $\{x_n\}_{n=1}^\infty$ and a sequence
$\{p_n\}_{n=1}^\infty$ of mutually disjoint projections in
$\m$ such that:
$$\lim_{n \to \infty} \Vert z_n - p_nz_np_n \Vert=0.$$
Note that $\liminf_{n \to \infty} \Vert p_nz_np_n \Vert >0$. 
By taking a subsequence (if necessary), we will assume that for
every $n \geq 1$,
$$\frac{\Vert z_n -p_nz_np_n \Vert}{\Vert p_nz_np_n\Vert} \leq
\varepsilon_n 2^{-n}.$$
For $n\geq 1$, set $y_n:= z_n/\Vert p_nz_np_n \Vert$.
If $(a_n)_n$ is a finite sequence of scalars then
\begin{equation*}
\begin{split}
\left\Vert \sum_n a_ny_n \right\Vert &\leq
\sum_n \vert a_n \vert \cdot 
\left\Vert y_n - p_nz_np_n\right\Vert  +
 \left(\sum_n \vert a_n \vert^p\right)^{1/p} \cr
&\leq \left(\sum_n \vert a_n \vert^p \varepsilon_{n}^p\right)^{1/p}
\cdot \left(\sum_n  2^{-nq} \right)^{1/q} +
\left(\sum_n \vert a_n \vert^p\right)^{1/p}
\end{split}
\end{equation*}
where $1/p +1/q =1$. This shows that
$$\left\Vert \sum_n a_ny_n \right\Vert \leq 
\left(\sum_n \vert a_n \vert^p\right)^{1/p} +
\left(\sum_n \vert a_n \vert^p\varepsilon_{n}^p\right)^{1/p}.$$
The other inequality can be obtain with similar estimate.
\end{proof}

\bigskip

\medskip

\noindent
{\bf Acknowledgement.} The author wishes to thank  F. Lust-Piquard and
G. Pisier for insightful discussions conserning the non-commutative
Khintchine's inequalities.

%\bibliography{narciref}
%\bibliographystyle{amsplain}

\providecommand{\bysame}{\leavevmode\hbox to3em{\hrulefill}\thinspace}

\end{document}